\documentclass[11pt]{article}
\usepackage[utf8]{inputenc}
\usepackage[T1]{fontenc}
\usepackage{lmodern}
\usepackage[a4paper]{geometry}
\usepackage[english]{babel}
\usepackage{pifont}
\usepackage{extsizes}
\usepackage{xcolor}
\usepackage{fancybox}
\usepackage{graphicx}
\usepackage{sectsty}
\usepackage{enumerate}
\usepackage{multicol}
\usepackage[english]{varioref}
\usepackage[fleqn]{amsmath}
\usepackage{mathtools}
\usepackage{amssymb}
\usepackage{amsmath}
\usepackage{mathrsfs}
\usepackage{bm}
\usepackage{amsthm}
\usepackage{blkarray}
\usepackage{subfigure}
\usepackage[babel=true]{csquotes}
\theoremstyle{break}
\newtheorem{Th}{Theorem}

\theoremstyle{definition}

\newcommand{\conv}[2][n]{\underset{#1\rightarrow #2}{\longrightarrow}}
\newcommand{\eq}[2][n]{\underset{#1\rightarrow #2}{\sim}}

\newtheorem{Le}{Lemma}

\newcommand{\sgn}{\text{sgn}}

\newcommand{\RR}{\mathbb{R}}
\newcommand{\ZZ}{\mathbb{Z}}

\newcommand{\NN}{\mathbb{N}}

\newcommand{\ind}[1]{\mathbf{1}_{\{#1\}}\,}
\newcommand{\PPP}[1]{\mathbb{P}\left(#1 \right)}%???crit la probabilit??? automatiquement avec la parenth???se

\newenvironment{prooft}[1]{\vskip 2mm\noindent {\bf Proof of #1.}}
                    {\hfill $\square$ \vskip 2mm \noindent}
\date{} 

\begin{document}
\title{Extremes for transient random walks in random sceneries under weak independence conditions}
\author{Nicolas Chenavier\footnote{Universit\'e du Littoral C\^ote d'Opale, Laboratoire de Math\'ematiques Pures et Appliqu\'ees J. Liouville, EA2597. Mail: nicolas.chenavier@univ-littoral.fr}, Ahmad Darwiche \footnote{Universit\'e du Littoral C\^ote d'Opale, Laboratoire de Math\'ematiques Pures et Appliqu\'ees J. Liouville, EA2597. Mail: Ahmad.Darwiche@univ-littoral.fr}}
\maketitle
\begin{abstract}
  Let $\{\xi(k), k \in \mathbb{Z} \}$ be a stationary sequence of random variables with conditions of type $D(u_n)$ and $D'(u_n)$. Let $\{S_n, n \in \mathbb{N} \}$ be a transient random walk in the domain of attraction of a stable law. We provide a limit theorem for the maximum of the first $n$ terms of the sequence  $\{\xi(S_n), n \in \mathbb{N} \}$ as $n$ goes to infinity. This paper extends a result due to Franke and Saigo who dealt with the case where the sequence $\{\xi(k), k \in \mathbb{Z} \}$ is i.i.d.
\end{abstract}


\section{Introduction} 
In the 1940s, Extreme Value Theory has been developed in the context of independent and identically distributed (i.i.d) random variables by Gnedenko \cite{Gnedenko}. It is straightforward that if $\{\xi(k), k\in \mathbb{Z}\}$ is a sequence of i.i.d random variables then the following property holds: for any sequence of real numbers $(u_n)$, and for $\tau >0$,
\[n\PPP{\xi(0)>u_n}\conv[n]{\infty}\tau \Longrightarrow \PPP{\max_{k\leq n}\xi(k)\leq u_n}\conv[n]{\infty}e^{-\tau}.\] 
The above property has been extended for sequences of dependent random variables satisfying the so-called conditions $D(u_n)$ and $D'(u_n)$ of Leadbetter \cite{L1,L2}. 

More recently, Franke and Saigo \cite{franke_saigo2009bis, franke_saigo_2009} have investigated extremes for a sequence of dependent random variables which do not satisfy the conditions $D(u_n)$ and $D'(u_n)$. More precisely, they consider the following problem. Let $\{X_{k}, k \in {\NN_+} \}$ be a sequence of centered, integer-valued i.i.d random variables and let $S_{n}=X_{1}+\dots+X_{n}$, $n\in \NN_+$. Assume that $\{X_{k}, k \in \NN_+ \}$ is in the domain of attraction of a stable law, i.e. for each $x\in \RR$, 
\begin{equation*}
\PPP{n^{-\frac{1}{\alpha}}S_{n}\leq x}  \conv[n]{\infty}  F_{\alpha}(x),    
\end{equation*} where $F_{\alpha}$ is the distribution function of a stable law with characteristic function given by
\[\varphi(\theta)=\exp(-|\theta|^{\alpha}(C_{1}+iC_{2}\sgn \theta)), \ \alpha \in (0,2].\] The sequence $\{S_n, n\in \NN_+\}$ is referred to as a random walk. When $\alpha<1$ (resp. $\alpha>1$), it is known that this random walk is transient (resp. recurrent) \cite{KS,LeGRo}. Now, let $\{\xi(k), k \in \mathbb{Z} \}$ be a family of $\mathbb{R}$-valued i.i.d random variables independent of the sequence  $\{X_{k}, k \in \NN_+ \}$. In the sense of \cite{franke_saigo_2009}, the sequence $\{\xi(S_n), n \in \NN_+ \}$ is called a random walk in a random scenery. Franke and Saigo derive limit theorems for the maximum of the first $n$ terms of $\{\xi(S_n), n \in \NN_+ \}$ as $n$ goes to infinity. An adaptation of Theorem 1 in \cite{franke_saigo_2009} shows that in the transient case, i.e. $\alpha<1$, the following property holds: if $n\PPP{\xi(0)>u_n}\conv[n]{\infty}\tau$ for some sequence $(u_n)$ and for some $\tau>0$, then \begin{equation}\label{eq:frankesaigo}  \PPP{\max_{k\leq n}\xi(S_k)\leq u_n} \conv[n]{\infty}e^{-\tau q},
\end{equation} where $q=\PPP{\forall k\in \NN_+, S_k\neq 0}$. Notice that $q>0$ because the random walk $\{S_n, n\in \NN_+\}$ is transient.  According to a result due to Le Gall and Rosen \cite{LeGRo}, the number $q$ can be also expressed as
\begin{equation}\label{eq:range} q=\lim_{n\rightarrow \infty}\frac{R_n}{n} \quad \text{a.s.},
\end{equation} where $R_n=\#\{S_1,\ldots, S_n\}$ is the range of the random walk.  

In this paper, we extend \eqref{eq:frankesaigo} to sequences $\{\xi(k), k\in \ZZ\}$ which are not necessarily i.i.d. but which only satisfy conditions of type $D(u_n)$ and $D'(u_n)$. More precisely, we consider the following problem. Let $S_n=X_1+\cdots +X_n$, where $\{X_{k}, k \in \NN_+ \}$ is a sequence satisfying the same properties as above, i.e. a sequence of centered, integer-valued i.i.d random variables in the domain of attraction of a $\alpha$-stable law, with $\alpha<1$. Let $\{\xi(k), k\in \ZZ\}$ be a stationary sequence of random variables independent of $\{X_{k}, k \in \NN_+ \}$. Assume that there exist a sequence $(u_n)$ such that \begin{equation}\label{eq:deftau}  n\PPP{\xi>u_n}\conv[n]{\infty}\tau,
\end{equation}  for some $\tau>0$, {where $\xi$ has the same distribution as $\xi(k)$, $k\in \ZZ$.} In the following, the sequence $\{\xi(k), k\in \ZZ\}$ is supposed to satisfy conditions of type $D(u_n)$ and $D'(u_n)$. Roughly, the condition $D(u_n)$ (see e.g. p29 in \cite{lucarini})  is a weak mixing property for the tails of the joint distributions. To introduce it, we write  for each $i_1< \cdots < i_p$ and for each $u\in \RR$,  \[F_{i_{1},\ldots,i_{p}}(u) = \PPP{\xi(i_1)\leq u, \ldots, \xi(i_p)\leq u}.\]  

\paragraph{Condition $D(u_n)$}  We say that $\{\xi(k), k \in \ZZ\}$ satisfies the condition $D(u_{n})$ if there exist a sequence $(\alpha_{n,l})_{(n,l)\in \NN^2}$ and a sequence $(l_n)$ of positive integers such that $\alpha_{n,l_{n}}\rightarrow 0$ as $n$ goes to infinity, $l_{n}=o(n)$, and 
\begin{equation*} 
|F_{i_{1},\ldots,i_{p}, j_{1},\ldots, j_{p'}}(u_{n})-F_{i_{1},\ldots,i_{p}}(u_{n})F_{j_{1},\ldots, j_{p'}}(u_{n})|\leq \alpha_{n,l}
\end{equation*}
 for any integers $i_{1}<\dots <i_{p}<j_{1}< \dots<j_{p'}$ such that $j_{1}-i_{p}\geq l$. Notice that the bound holds uniformly in $p$ and $p'$.

The condition $D'(u_n)$ (see e.g. p29 in \cite{lucarini}) is a local type property and precludes the existence of clusters of exceedances. To introduce it, we consider a sequence $(k_n)$ such that
\begin{equation}
\label{eq:defrn}
k_n\conv[n]{\infty}\infty, \quad  \frac{n^2}{k_n} \alpha_{n,l_{n}}\conv[n]{\infty}0, \quad k_nl_n=o(n),
\end{equation}
where $(l_n)$ and $(\alpha_{n,l})_{(n,l)\in \NN^2}$ are the same as in condition $D(u_n)$.

\paragraph{Condition $D'(u_n)$}
We say that $\{\xi(k), k \in \ZZ\}$ satisfies the condition $D'(u_{n})$  if there exists a sequence of integers $(k_n)$ satisfying \eqref{eq:defrn} such that 
\[\lim\limits_{n\rightarrow \infty} n \sum_{j=1}^{\lfloor n/k_n\rfloor}\PPP{\xi(0)>u_{n},\xi(j)>u_n}=0.\]
In the classical literature, the sequences $(\alpha_{n,l})_{(n,l)\in \NN^2}$ and $(k_n)$ only satisfy $k_n\alpha_{n,l_n}\conv[n]{\infty}0$ (see e.g. (3.2.1) in \cite{lucarini}) whereas in \eqref{eq:defrn} we have assumed that  $\frac{n^2}{k_n} \alpha_{n,l_{n}}\conv[n]{\infty}0$. In this sense, the condition $D'(u_n)$ as written above is slightly more restrictive than the usual condition $D'(u_n)$.  We are now prepared to state our main theorem. 

\begin{Th}
\label{Th:maintheorem}
Let $\{S_n, n \in \NN_+\}$ be as above and let $\{\xi(k), k \in \ZZ\}$ be a stationary sequence of random variables such that $n\PPP{\xi(0)>u_n}\conv[n]{\infty}\tau$, for some sequence $(u_n)$ and $\tau \geq 0$. Assume that the conditions $D(u_{n})$ and $D'(u_{n})$ hold. Then for almost all realization of $\{S_n, n \in \NN_+\}$, 
\[\PPP{\max_{k\leq n}\xi(S_k)\leq u_{n}} \conv[n]{\infty}e^{-\tau q},\]
where $q=\PPP{\forall k\in \NN_+, S_k\neq 0}$.
\end{Th}
{The above result extends \eqref{eq:frankesaigo} to} sequences $\{\xi(k),k\in \ZZ\}$ which only satisfy the conditions $D(u_n)$ and $D'(u_n)$. Notice that our theorem is stated when $\alpha<1$, but it remains true when $\alpha\geq 1$. We did not deal with this case because, when the random walk is recurrent, the number $q$ equals 0 and the limit is degenerate, i.e.  $\PPP{\max_{k\leq n}\xi(S_k)\leq u_{n}} \conv[n]{\infty}1$. 

The main idea to derive Theorem \ref{Th:maintheorem} is to adapt \cite{L2} to our context. We think that our method combined with Kallenberg's theorem ensures that the point process of exceedances converges to a Poisson point process, in the same spirit as Theorem 3 in \cite{franke_saigo_2009}. More precisely, if the threshold is of the form $u_n=u_n(x)=a_nx+b_n$, for some $x\in \RR$, and if we let $\tau_k=\inf\{m\in \NN_+, \#\{S_1, \ldots, S_m\}\geq k\}$, then the point process \[\Phi_n=\left\{\left(\frac{\tau_k}{n}, \frac{\xi(S_{\tau_k}) - b_{\lfloor qn\rfloor}}{a_{\lfloor qn\rfloor}}\right), k\geq 1\right\}\] converges to a Poisson point process with explicit intensity measure. 

The rest of the paper is devoted to the proof of Theorem \ref{Th:maintheorem}.
%%%%%%%%%%%%%%%%%%%%%%%%%%%%%%%%%%%%%%%%%%%%%%%%%%%%%%%%%%%%%%%%%%%%%%%%%%%%%%%%

\section{Proof of Theorem \ref{Th:maintheorem}}

 The main idea is to adapt \cite{L2} to our context. To do it, let $(k_n)$, $(l_n)$ be as in  \eqref{eq:defrn}. {For $n$ large enough, let $r_n=\lfloor\frac{n}{k_n-1}\rfloor+1$.}  Given a realization $\{S_n, n\in \NN_+\}$ of the random walk, we write $\mathcal{S}_{n}=\{S_{1},\ldots,S_{n}\}$ and $R_n=\#\mathcal{S}_n$.  To capture the fact that the random scenery $\{\xi(k),  k\in \ZZ\}$ satisfies the condition $D(u_n)$, we construct blocks and stripes as follows. 
Let \begin{equation} \label{eq:defKn} K_n=\left\lfloor\frac{R_n}{r_n}\right\rfloor +1.\end{equation}
There exists a unique $K_n$-tuple of subsets $B_i\subset \mathcal{S}_n$, $i\leq K_n$, such that the following properties hold: $\bigcup_{j\leq K_n}B_j=\mathcal{S}_n$, $\#B_i=r_n$ and $\max B_i<\min B_{i+1}$ for all $i\leq K_n-1$. Notice that { $K_n\leq k_n$ and $\#B_{K_n}=R_n-(K_n-1)\cdot r_n$ almost surely (a.s.)}. The sets $B_{j}$, $j\leq K_n$, are referred to as \textit{blocks}. For each $j\leq K_n$, we also denote by $L_j$ the family consisting of the $l_n$ largest  terms of $B_j$ (e.g. if $B_j=\{x_1,\ldots, x_{r_n}\}$, with $x_1<\cdots <x_{r_n}$, $j\leq K_n-1$, then $L_j=\{x_{r_n-l_n+1}, \ldots, x_{r_n}\}$). When $j=K_n$, we take the convention $L_{K_n}=\emptyset$ if $\#B_{K_n}<l_n$.  The set $L_{j}$  is referred to as a \textit{stripe}, {and the union of the stripes is denoted $\mathcal{L}_n=\bigcup_{j\leq K_n}L_{j}$. In the rest of the paper, we write $M_B=\max_{k\in B}\xi(k)$ for all subset $B\subset \ZZ$.} To prove Theorem \ref{Th:maintheorem}, we will use the following lemma.

\begin{Le} \label{Le:stripes}
With the above notation, we have for almost all realization of $\{S_n, n\in \NN_+\}$,
\begin{enumerate}[(i)]
\item $	\PPP{M_{\mathcal{S}_{n}}\leq u_{n}} - \PPP{M_{\mathcal{S}_{n}\setminus \mathcal{L}_{n}}\leq u_{n}}\conv[n]{\infty}0$;
\item $\PPP{M_{\mathcal{S}_{n}\setminus \mathcal{L}_{n}}\leq u_{n}}-\prod_{j\leq K_{n}} \PPP{M_{B_{j}\setminus \mathcal{L}_{n}}\leq u_{n}}\conv[n]{\infty}0$;
\item  $ \prod_{j\leq K_{n}} \PPP{M_{B_{j}\setminus \mathcal{L}_{n}}\leq u_{n}}- \prod_{j\leq K_{n}} \PPP{M_{B_{j}}\leq u_{n}} \conv[n]{\infty}0$.
\end{enumerate}
\end{Le}
The first and the third assertions mean that, asymptotically, the maximum is not affected if we remove the sites which belong to one of the stripes. Roughly, this comes from the fact that the size of the stripes is negligible compared to the size of the blocks. The second assertion is a consequence of the fact that the sequence $\{\xi(k), k\in \ZZ\}$ satisfies the condition $D(u_n)$. To derive Theorem \ref{Th:maintheorem}, we will also use the following lemma.
\begin{Le}\label{Le:cluster}
With the above notation, we have for almost all realization of $\{S_n, n\in \NN_+\}$, 
\[ \prod_{j\leq K_{n}} \PPP{M_{B_{j}}\leq u_{n}}\conv[n]{\infty} e^{-\tau q}.\]
\end{Le}
\begin{prooft}{Lemma \ref{Le:stripes}}
First we prove (i). To do it, for all $n\in \NN_+$, we write
\begin{align}
\label{eq:majstripe}
0\leq   \PPP{M_{\mathcal{S}_{n}\setminus \mathcal{L}_{n}}\leq u_{n}}-\PPP{ M_{\mathcal{S}_{n}}\leq u_{n}} & \leq \PPP{M_{ \mathcal{L}_{n} }> u_{n}} \notag \\
  &\leq \# \mathcal{L}_{n} \PPP{\xi > u_{n}} \notag \\
  &\leq K_nl_n\PPP{\xi > u_{n}}\notag\\
  & \leq k_n l_n \PPP{\xi > u_{n}}.
\end{align}
Since $k_n l_n = o(n)$ and $n\PPP{\xi > u_{n}}\conv[n]{\infty} \tau$, we have $k_n l_n\PPP{\xi > u_{n}}\conv[n]{\infty}0$. This together with \eqref{eq:majstripe} concludes the proof of (i).

Now we  prove (ii). Noticing that $\{M_{\mathcal{S}_{n}\setminus \mathcal{L}_{n}}\leq u_{n}\}=\bigcap_{j\leq K_n}\{M_{B_{j}\setminus \mathcal{L}_{n}}\leq u_{n}\}$ and bounding $\PPP{M_{ B_{K_n}\setminus \mathcal{L}_{n}}\leq u_{n}}$ by $1$, we have 
\begin{multline*}
  \left| \PPP{M_{\mathcal{S}_{n}\setminus \mathcal{L}_{n}}\leq u_{n}}-\prod_{j\leq K_{n}} \PPP{M_{B_{j}\setminus \mathcal{L}_{n}}\leq u_{n}}\right|\\ 
  \leq \left|\PPP{\bigcap_{j\leq K_n}M_{B_{j}\setminus \mathcal{L}_{n}}\leq u_{n}}- \PPP{\bigcap_{j\leq K_n-1}M_{B_{j}\setminus \mathcal{L}_{n}}\leq u_{n}}\PPP{M_{ B_{K_n}\setminus \mathcal{L}_{n}}\leq u_{n}}\right|\\
 +\left| \PPP{\bigcap_{j\leq K_n-1}M_{B_{j}\setminus \mathcal{L}_{n}}\leq u_{n}}-\prod_{j\leq K_{n}-1} \PPP{M_{B_{j}\setminus \mathcal{L}_{n}}\leq u_{n}}\right|.
\end{multline*}
{It follows from the definition of $\mathcal{L}_n$ that $\inf_{j\leq K_n-1}d(B_{K_n}\setminus \mathcal{L}_n, B_{j}\setminus \mathcal{L}_n)\geq l_n$, where $d(A,B)$ denotes the distance between any pairs of sets $A,B\subset \RR$.} Thanks to the condition $D(u_{n})$, {this gives}
\begin{multline*}
 \left| \PPP{M_{\mathcal{S}_{n}\setminus \mathcal{L}_{n}}\leq u_{n}}-\prod_{j\leq K_{n}} \PPP{M_{B_{j}\setminus \mathcal{L}_{n}}\leq u_{n}}\right|\\
 \leq \alpha_{n,l_{n}}+\left|\PPP{\bigcap_{j\leq K_n-1}M_{B_{j}\setminus \mathcal{L}_{n}}\leq u_{n}}-\prod_{j\leq K_{n}-1} \PPP{M_{B_{j}\setminus \mathcal{L}_{n}}\leq u_{n}}\right|.
\end{multline*}
By induction, we have
\begin{align*}
  \left| \PPP{M_{\mathcal{S}_{n}\setminus \mathcal{L}_{n}}\leq u_{n}}-\prod_{j\leq K_{n}} \PPP{M_{B_{j}\setminus \mathcal{L}_{n}}\leq u_{n}}\right| \leq (K_{n}-1) \alpha_{n,l_{n}}\leq k_n\alpha_{n,l_{n}},
\end{align*}
which converges to 0 as $n$ goes to infinity. This concludes the proof of (ii). 

It remains to prove (iii). To do it, notice that, for $n$ large enough, $\PPP{M_{B_{j}}\leq u_{n}}\neq 0$ {because}  $\PPP{M_{B_{j}}>u_{n}} \leq r_{n} \PPP{\xi>u_{n}} \conv[n]{\infty}0$. This allows us to write for $n$ large enough,
\begin{multline}
\label{eq:majiii}
\left| \prod_{j\leq K_{n}} \PPP{M_{B_{j}\setminus \mathcal{L}_{n}}\leq u_{n}}- \prod_{j\leq K_{n}} \PPP{M_{B_{j}}\leq u_{n}} \right| \\
\begin{split}
&=\prod_{j\leq K_{n}} \PPP{M_{B_{j}}\leq u_{n}} \left(\prod_{j\leq K_{n}} \frac{\PPP{M_{B_{j}\setminus \mathcal{L}_{n}}\leq u_{n}}}{\PPP{M_{B_{j}}\leq u_{n}}}-1\right),\\
&\leq \prod_{j\leq K_{n}} \frac{\PPP{M_{B_{j}\setminus \mathcal{L}_{n}}\leq u_{n}}}{\PPP{M_{B_{j}}\leq u_{n}}}-1. 
\end{split}
\end{multline}
Now, let $j\leq K_n$ be fixed. Because $\#\{B_j\cap \mathcal{L}_n\}\leq l_n$, we can prove, in the same spirit as  \eqref{eq:majstripe}, that  
 \[\PPP{M_{B_{j}\setminus \mathcal{L}_{n}}\leq u_{n}}-\PPP{M_{B_{j}}\leq u_{n}}\leq l_{n}\PPP{\xi>u_{n}}.\] Adapting \eqref{eq:majstripe} again, we also have \begin{equation}\label{eq:majmaxBj} \PPP{M_{B_{j}}>u_{n}} \leq r_{n} \PPP{\xi>u_{n}}.\end{equation} This implies that 
 \begin{align*}
  \frac{\PPP{M_{B_{j}\setminus \mathcal{L}_{n}}\leq u_{n}}}{\PPP{M_{B_{j}}\leq u_{n}}}
& = 1 + \frac{\PPP{M_{B_{j}\setminus \mathcal{L}_{n}}\leq u_{n}}- \PPP{M_{B_{j}}\leq u_{n}}}{1-\PPP{M_{B_j}>u_n}}  \\
& \leq 1+ \frac{l_{n}\PPP{\xi>u_{n}}}{1-r_{n} \PPP{\xi>u_{n}}}\\
& = 1+ \frac{\tau l_{n}}{n}+o\left(\frac{l_n}{n}\right),  
\end{align*}
where the last line comes from the fact that $\PPP{\xi>u_n}\eq[n]{\infty}\frac{\tau}{n}$. This together with \eqref{eq:majiii} and the fact that $K_n\leq k_n$ implies
\begin{align*}
\left| \prod_{j\leq K_{n}} \PPP{M_{B_{j}\setminus \mathcal{L}_{n}}\leq u_{n}}- \prod_{j\leq K_{n}} \PPP{M_{B_{j}}\leq u_{n}} \right| & \leq  \left( 1+\frac{\tau l_n}{n} +o\left(\frac{l_n}{n}\right) \right)^{k_n}-1\\
& \eq[n]{\infty}\frac{\tau k_nl_n}{n}. 
\end{align*}
The last term converges to 0 as $n$ goes to infinity since $k_nl_n=o(n)$. This concludes the proof of (iii). 
\end{prooft}
%%%%%%%%%%%%%%%%%%%%%%%%%%%%%%%%%%%%%%%%%%%%%%%%%%%%%%%%%%%%%%%%%%%%%%%%%%%%%
\begin{prooft}{Lemma \ref{Le:cluster}}
First, we provide a lower-bound for  $\prod_{j\leq K_{n}} \PPP{M_{B_{j}}\leq u_{n}}$. To do it, {for all $n\in \NN_+$,} we write
\begin{align*} \prod_{j\leq K_{n}} \PPP{M_{B_{j}}\leq u_{n}} &  =\exp\left(\sum_{j \leq K_{n}}\log\left(1-\PPP{M_{B_{j}}>u_{n}}\right)\right)\\
& \geq \exp\left(K_n\log\left(1- r_{n} \PPP{\xi>u_{n}}\right)\right),\end{align*}
where the last line comes from \eqref{eq:majmaxBj}. To deal with the  right-hand side, we notice that \[\log\left(1- r_{n} \PPP{\xi>u_{n}}\right) \eq[n]{\infty}-r_n\PPP{\xi>u_{n}}\] since $r_{n} \PPP{\xi>u_{n}}$ converges to 0. It follows from  \eqref{eq:range}, \eqref{eq:deftau} and  \eqref{eq:defKn} that, for almost all realization of $\{S_n, n\in \NN_+\}$,  \[\lim_{n\rightarrow\infty}\exp\left(K_n\log\left(1- r_{n} \PPP{\xi>u_{n}}\right)\right) = \exp(-\tau q).\]  Therefore a.s.
\[ \liminf_{n\rightarrow\infty}\prod_{j\leq K_{n}} \PPP{M_{B_{j}}\leq u_{n}}\geq  \exp(-\tau q).\]

Now, we provide an upper-bound for $\prod_{j\leq K_{n}} \PPP{M_{B_{j}}\leq u_{n}}$. To do it, we write
\begin{align*} \prod_{j\leq K_{n}} \PPP{M_{B_{j}}\leq u_{n}} &  =\exp\left(\sum_{j \leq K_{n}}\log\left(1-\PPP{M_{B_{j}}>u_{n}}\right)\right)\\
& \leq \exp\left(-\sum_{j\leq K_n}\PPP{M_{B_j}>u_n}\right).
\end{align*}
This together with the Bonferroni inequalities (see e.g. p110 in \cite{Fel}), implies that 
\begin{multline*}
 \prod_{j\leq K_{n}} \PPP{M_{B_{j}}\leq u_{n}} \\
 \leq \exp\left( -{(K_n-1)}r_n \PPP{\xi>u_{n}}{+}\sum_{j\leq K_n} \sum_{\alpha < \beta; \alpha, \beta \in B_{j}}\PPP{\xi(\alpha)>u_{n},\xi(\beta)>u_{n}} \right).
\end{multline*} 
Since $K_nr_n\PPP{\xi>u_n}\conv[n]{\infty}\tau q$ a.s., we have a.s. 
\[\limsup_{n\rightarrow\infty} \prod_{j\leq K_{n}} \PPP{M_{B_{j}}\leq u_{n}}\leq  \exp\left( -\tau q{+\limsup_{n\rightarrow\infty}}\sum_{j\leq K_n} \sum_{\alpha < \beta; \alpha, \beta \in B_{j}}\PPP{\xi(\alpha)>u_{n},\xi(\beta)>u_{n}} \right).  \] Therefore, it is enough to prove that
\[\sum_{j\leq K_n} \sum_{\alpha < \beta; \alpha, \beta \in B_{j}}\PPP{\xi(\alpha)>u_{n},\xi(\beta)>u_{n}} \conv[n]{\infty}0.\]
To do it, we write the sum appearing in the above equation into two terms: the first one deals with the case when 
$\beta-\alpha<r_n$, and the second one deals with the opposite. For the first term, we {use the fact that the sequence $\{\xi(k), k\in \ZZ\}$ is stationary. This gives}  
\begin{align*}
\sum_{j\leq K_n} \sum_{\alpha < \beta; \alpha, \beta \in B_{j}}\PPP{\xi(\alpha)>u_{n},\xi(\beta)>u_{n}}&\ind{\beta-\alpha<r_n}\\
&\leq  \sum_{j\leq K_n}\sum_{\beta \in B_{j}} \sum_{k=1}^{r_{n}} \PPP{\xi(0)>u_{n},\xi(k)>u_{n}} \\
&\leq k_nr_n \sum_{k=1}^{r_{n}} \PPP{\xi(0)>u_{n},\xi(k)>u_{n}} 
\end{align*} 
The last term converges to 0 as $n$ goes to infinity according to the condition $D'(u_n)$ and the fact that $k_nr_n\eq[n]{\infty}n$.  Now, we deal with the same series but this time by replacing  $\ind{\beta-\alpha<r_n}$ by $\ind{\beta-\alpha\geq r_n}$. We have
\begin{multline}
\label{eq:majbetasupalpha}
  \sum_{j\leq K_n} \sum_{\alpha < \beta; \alpha, \beta \in B_j} \PPP{\xi(\alpha)>u_n,\xi(\beta)>u_n}\ind{\beta-\alpha\geq r_n}\leq \sum_{j\leq K_n}  \sum_{\alpha < \beta; \alpha,\beta \in B_j} \PPP{\xi>u_n}^2\\
    +    \sum_{j\leq K_n} \sum_  {\alpha < \beta; \alpha,\beta \in B_{j}}\left|\PPP{\xi(\alpha)>u_{n},\xi(\beta)>u_{n}}-\PPP{\xi>u_n}^2\right|\ind{\beta-\alpha\geq r_n}.
\end{multline}
We prove below that the two terms of the right-hand side converge to 0. For the first one, we {have}
\[\sum_{j\leq K_n}  \sum_{\alpha < \beta; \alpha,\beta \in B_j} \PPP{\xi>u_n}^2\leq k_nr_n^2\PPP{\xi>u_n}^2\eq[n]{\infty}\tau^2\frac{r_n}{n}.\]
The last term converges to 0 according to \eqref{eq:defrn}. To deal with the second term of \eqref{eq:majbetasupalpha}, we use the condition $D(u_n)$. This gives
\begin{multline*}  \sum_{j\leq K_n} \sum_  {\alpha < \beta; \alpha,\beta \in B_{j}}\left|\PPP{\xi(\alpha)>u_{n},\xi(\beta)>u_{n}}-\PPP{\xi>u_n}^2\right|\ind{\beta-\alpha\geq r_n}   \\ 
\begin{split}  
& \leq  \sum_{j\leq K_n} \sum_  {\alpha < \beta; \alpha,\beta \in B_{j}}\left|\PPP{\xi(\alpha)>u_{n},\xi(\beta)>u_{n}}-\PPP{\xi>u_n}^2\right|\ind{\beta-\alpha\geq l_n}\\
& \leq k_nr_n^2\alpha_{n,l_n}
\eq[n]{\infty}\frac{n^2}{k_n}\alpha_{n,l_n},
\end{split}
\end{multline*}
which converges to 0 as $n$ goes to infinity according to \eqref{eq:defrn}. This concludes the proof of Lemma \ref{Le:cluster}.
\end{prooft}
Theorem \ref{Th:maintheorem} follows directly from Lemmas \ref{Le:stripes} and \ref{Le:cluster}.

%\bibliographystyle{abbrv}
%\bibliography{BiblioVE}

\end{document}